\documentclass[10pt, a4paper]{article}
\usepackage{latexsym}
\usepackage{indentfirst}
\usepackage{amsmath}
\usepackage{amssymb}

\begin{document}

\title{On the computation of the term $w_{21}z^2\overline{z}$ of the series defining the
center manifold for a scalar delay differential equation}
\author{Anca-Veronica Ion\\
"Gh. Mihoc-C. Iacob" Institute of Mathematical Statistics\\
and Applied Mathematics of the Romanian Academy, \\
13, Calea 13 Septembrie, Bucharest, Romania \\
anca.veronica.ion@ima.ro}

\date{}
\maketitle

\begin{abstract}
In computing the third order terms of the series of powers of the
center manifold at an equilibrium point of a scalar delay
differential equation, some problems occur at the term
$w_{21}z^2\overline{z}.$ More precisely, in order to determine the
values at 0, respectively $-r$ of the function $w_{21}(\,.\,),$ an
algebraic system of equations must be solved. We show that the two
equations are dependent, hence the system has an infinity of
solutions. Then we show how we can overcome this lack of
uniqueness and provide a formula for $w_{21}(0).$
\end{abstract}

\section{Introduction}

We consider the delay differential equation
\begin{equation}\label{eq}\dot{x}(t)=f(x(t),x(t-r)),
\end{equation}
where the function $f$ is a scalar function, $f(0,0)=0,$ and $r>0.$
For the sake of simplicity we assume $f\in C^{\infty}(D)$ ($(0,0)\in
D\subset\mathbb{R}^2,$ $D$ open) but the results below are valid
also for a lower class of smoothness (e.g. $C^k,\,k>3$).

As usual in the study of delay differential equations \cite{Ha},
\cite{HaL} we consider the function space
$\mathcal{B}=\left\{\psi:[-r,0]\mapsto \mathbb{R},\, \psi\,
\mathrm{is\, continuous\, on\,}[-r,0]\, \right\}$ and its
complexification, $\mathcal{B}_C=\mathcal{B}+i\mathcal{B}$. For a
continuous function $x:[-r,T]\mapsto \mathbb{R},\;T>0,$ we denote by
$x_t$ the element of $\mathcal{B}$ defined by
$x_t(s)=x(t+s),\;s\in[-r,0].$

We also assume that the function $f$ is such that, for every
$\phi\in \mathcal{B}$ the above equation, with the initial condition
$x_0=\phi$ has an unique solution $x(t,\phi),\;t\geq 0$ and, thus,
the equation generates a semi-dynamical system $\{T(t)\}_{t\geq 0}$
on $\mathcal{B},$ by $T(t)\phi=x_t(\phi).$

In order to simplify the sequel developments, we write eq. \eqref{eq} in the form
\begin{equation}\label{s-eq}\dot{x}(t)=Ax(t)+Bx(t-r)+\widehat{f}(x(t),\,x(t-r)),
\end{equation}
 $\widehat{f}$ being the nonlinear part of $f$,
\[\widehat{f}(y,z)=\sum_{j,k\geq0,
j+k\geq2}\frac{1}{j!k!}C_{j,k}y^jz^k.
\]

 The linearized equation attached to the previous, that is
\begin{equation}\label{lin-eq}\dot{x}(t)=Ax(t)+Bx(t-r),
\end{equation}
generates a semigroup of operators on $\mathcal{B}.$ The eigenvalues of the infinitesimal generator of this semigroup are the solutions of
\begin{equation}\label{ch-eq}\lambda-A-e^{-\lambda r}B=0.
\end{equation}
\textit{We assume that \eqref{ch-eq} has a pair of pure imaginary
complex conjugated solutions $\lambda_{1,2}=\pm \omega i,$ with
$\omega>0,$ and all other eigenvalues have negative real part. }

Consider the subspace of $\mathcal{B}_C$ spanned by $\varphi_{1,2}$ and denote it by $\mathcal{M}.$
The space $\mathcal{B}_C$ is decomposed as a direct sum $\mathcal{B}_C=\mathcal{M}\bigoplus \mathcal{N}$ with the help of a projector. This projector is built with the help of a bilinear form, and this one is constructed by using the ``adjoint'' problem \cite{Ha}, \cite{HaL}.

In the above conditions, it is proved \cite{HaL} that there is a
local invariant manifold, called the \textit{center manifold} that
is a smooth manifold, tangent to the space $\mathcal{M}$ at the
point $x=0$
 and it is the graph of a function
$w(\cdot)$ defined on a neighborhood of zero in
$\mathcal{M}$ and taking values in $\mathcal{N}$.
A point on the local invariant manifold has the form
$u\,\varphi_1+\overline{u}\,\varphi_2+w(u\,\varphi_1+\overline{u}\,\varphi_2).$

In order to approximate the function $w$ that defines the center manifold, we define $\widehat{w}(z,\overline{z}):=w(z\varphi_1+\overline{z}\varphi_2).$ For simplicity, we will drop the hat, and we will denote also by $w$ the new two-variables function defined above.

We write
\begin{equation}\label{cm}  w(z,\overline{z})=\sum_{i+j\geq 2}\frac{1}{i!j!}w_{i,j}z^i\overline{z}^j,\end{equation} where  $w_{i,j}\in \mathcal{B}_C$.
From the invariance of the center manifold, a set of differential
equations for the functions $w_{i,j}$ is obtained as well as a set
of conditions for the determination of the integration constants.
More precisely, these latter conditions have the form of linear
relations between $w_{i,j}(0)$ and $w_{i,j}(-r)$.

For the second order terms of the series, the computation does not
present any problems. This happens also for the third order terms,
excepting $w_{2,1}$.  For the function $w_{2,1}$ the determinant of
the linear algebraic system in $w_{i,j}(0)$ and $w_{i,j}(-r)$ is
zero - Proposition 3.1. We prove that the system has (an infinity
of) solutions - Proposition 3.2.

The new problem is how to select the proper values for $w_{2,1}(0)$
and $w_{2,1}(-r)$.

In Section 4 we consider a perturbed problem, depending on a small
parameter $\epsilon>0,$ such that, when $\epsilon \rightarrow 0,$
the unperturbed problem is obtained. We show that for this perturbed
problem the system of algebraic equations for the corresponding
$w_{\epsilon2,1}(0),$ $w_{\epsilon2,1}(-r)$ has unique solution and
that $w_{2,1}(0),$ $w_{2,1}(-r)$ can be obtained by taking the limit
when $\epsilon\rightarrow 0$ - Proposition 4.1 and Proposition 4.2.
A consequence of the proof is the fact that the limit does not
depend on the particular perturbed problem considered. A formula for
$w_{2,1}(0)$ is given in Proposition 4.3.

\section{The reduction of the problem to the center manifold}

For the presentation of our results we need the following framework
of \cite{Ha}, \cite{HaL}, \cite{Far}. More precisely, we intend to describe the method for
constructing the restriction of the problem to the center manifold.

For this, we also need  the space
$$\mathcal{B}_0=\left\{\psi:[-r,0]\mapsto \mathbb{R},\, \psi\,
\mathrm{is\, continuous\, on\,}[-r,0)\wedge\,\exists
\lim_{s\rightarrow 0}\psi(s)\in \mathbb{R}  \right\},$$ and its
complexification, $\mathcal{B}_{0C}.$ Notice that $\mathcal{B}_0$
consists of functions of the form $\psi=\varphi+\sigma d_{0},$ where
$\varphi\in \mathcal{B},$ $\sigma\in \mathbb{R}$ and
$d_0:[-r,0]\mapsto \mathbb{R},$
$$d_0(s)=\left\{\begin{array}{cc}
                0, & s\in[-r,0), \\
                1, & s=0,\\
              \end{array}\right.
$$
with norm given by $\|\psi\|=|\varphi|_0+|\sigma|.$

The linear part of the RHS of eq. \eqref{s-eq} can be written with
the help of linear operator $L:\mathcal{B}\mapsto \mathbb{R},$ given
by $L(\varphi)=A\varphi(0)+B\varphi(-r).$ For further use, we write
this operator, with the help of a Stieljes integral as
\[\quad L\varphi =\int_{-r}^0\varphi(\theta)d\eta(\theta),
\]
where
\begin{equation}\label{eta}\eta(s)=\left\{\begin{array}{c}
              -B,\quad s=-r; \\
                \quad\quad\;\, 0,\quad\; s\in (-r,0); \\
              A,\quad s=0.\\
            \end{array}\right.
\end{equation}

In \cite{Far}, the linear operator
$\mathcal{A}:C^{1}([-r,0],\mathbb{R})\subset\mathcal{B}_0\mapsto
\mathcal{B}_0$
\begin{equation}\mathcal{A}(\varphi)=\dot{\varphi}+d_0[L(\varphi)-\dot{\varphi}(0)]
\end{equation}
is defined and it is proved that this is the infinitesimal generator
of the semigroup of operators $\{S(t)\}_{t\geq 0}$ given by
$S(t)(\phi)=x_t(\phi),$ where $x(t,\phi)$ is the solution of
equation (\ref{lin-eq}) with the initial condition $x_0=\phi.$ Then
the nonlinear equation may be written as an equation in
$\mathcal{B}_0,$ that is
\begin{equation}\label{nonlineq3}\frac{dx_t}{dt}=\mathcal{A}(x_t)+d_0\widetilde{f}(x_t),
\end{equation}
where
$\widetilde{f}(\varphi)=\widehat{f}(\varphi(0),\,\varphi(-r)).$

In \cite{Ha}, \cite{HaL}, the adjoint equation (associated to the
linear equation \eqref{lin-eq}) is defined as
\[\dot{y}(s)=-Ay(s)-By(s+r).
\]
The corresponding characteristic equation is
\[\lambda+A+Be^{\lambda r}=0,
\]
and it is obvious that, together with eq. \eqref{ch-eq}, it admits the solutions $\pm \omega i.$
 The corresponding eigenfunctions  are $\psi_{1}(\zeta)=e^{-\omega i\zeta},\;\psi_{2}(\zeta)=e^{\omega i\zeta},\;\zeta\in[0,r].$

Also in \cite{Ha}, \cite{HaL} in order to construct a projector on
$\mathcal{M},$ the following bilinear form is defined on
$C([0,r],\mathbb{C})\times \mathcal{B}_C$
\begin{equation}\label{bf}\langle\psi,\varphi\rangle=\psi(0)\varphi(0)-\int_{-r}^0\int_{0}^\theta
\psi(\zeta-\theta)\varphi(\zeta)d\zeta d\eta(\theta)=
\end{equation}
\[=\psi(0)\varphi(0)+
B\int_{-r}^0\psi(\zeta+r)\varphi(\zeta)d\zeta,
\]
($\eta$ being the function defined in \eqref{eta}). Then linear
combinations of the functions $\psi_j, \;j=1,2,$ denoted by
$\Psi_i,\,i=1,2$ are constructed such that $\langle \Psi_i,
\varphi_j \rangle=\delta_{ij}.$ For this we determine the $2\times
2$ matrix $E,$ with elements
$e_{ij}=\langle\psi_i,\,\varphi_j\rangle$:

\[e_{11}=\langle \psi_1,\varphi_1\rangle=\psi_1(0)\varphi_1(0)+
B\int_{-r}^0e^{-i\omega(\theta+r)}e^{i\omega\theta}d\theta=1-(A-i\omega)r,
\]

\[e_{12}=\langle\psi_1,\varphi_2\rangle=\psi_1(0)\varphi_2(0)+
B\int_{-r}^0e^{-i\omega(\theta+r)}e^{-i\omega\theta}d\theta=0,
\]

\[e_{21}=\langle\psi_2,\varphi_1\rangle=\psi_2(0)\varphi_1(0)+
B\int_{-r}^0e^{i\omega(\theta+r)}e^{i\omega\theta}d\theta=0,
\]
\[e_{22}=\langle\psi_2,\varphi_2\rangle=\psi_2(0)\varphi_2(0)+
B\int_{-r}^0e^{i\omega(\theta+r)}e^{-i\omega\theta}d\theta=1-(A+i\omega)r.
\]

Then \[\left(%
\begin{array}{c}
  \Psi_1 \\
  \Psi_2 \\
\end{array}%
\right)=E^{-1}\left(%
\begin{array}{c}
  \psi_1 \\
  \psi_2 \\
\end{array}%
\right),
\]
and we obtain \begin{equation}\label{Psi}\Psi_1(\zeta)=\frac{e_{22}}{\det
E}\psi_1(\zeta)=\frac{1-(A+i\omega)r}{(1-Ar)^2+\omega^2r^2}e^{-\omega i \zeta},\quad
\Psi_2=\overline{\Psi}_1.
\end{equation}

The projector defined in \cite{Far} on $\mathcal{B}_{0C}$ and with
values in
$\mathcal{M}$  is given, for\\
$\psi=\phi+d_0\sigma\in\mathcal{B}_{0C}$ by
\begin{equation}\mathcal{P}(\psi)=\left(\langle\Psi_1,\phi\rangle+\Psi_1(0)\sigma\right)\varphi_1
+\left(\langle\Psi_2,\phi\rangle+\Psi_2(0)\sigma\right)\varphi_2.\end{equation}
If $\psi=\phi\in\mathcal{B}_C,$ ($\sigma=0$) we have
$\mathcal{P}(\phi)=\langle\Psi_1,\phi\rangle\varphi_1
+\langle\Psi_2,\phi\rangle\varphi_2.$

Now, for the solution $z(\cdot,\phi)$ of (\ref{nonlineq3}) with
initial condition $x_0=\phi,$ we can write
\[z_t=\varphi_1u_1(t)+\varphi_2u_2(t)+\mathbf{v}(t),\,
\]
with $u_1(t)=\langle\Psi_1,\,z_t
\rangle,\,\,u_2(t)=\langle\Psi_2,\,z_t
\rangle,\,\,\mathbf{v}=(I-\mathcal{P})z_t.$

In \cite{Far}  the equation (\ref{nonlineq3}) is projected by $\mathcal{P},$ to
obtain,
\begin{equation}\label{proj-eq}
\frac{d\mathbf{u}}{dt}=\Lambda\mathbf{u}+\Psi(0)\widetilde{f}(\varphi_1u_1+\varphi_2u_2+\mathbf{v}),
\end{equation}
where $$\mathbf{u}=\left(%
\begin{array}{c}
  u_1 \\
  u_2 \\
\end{array}%
\right), \,\,\Lambda=\left(%
\begin{array}{cc}
  \lambda_1 & 0 \\
  0 & \lambda_2 \\
\end{array}\right),\,\,\Psi=\left(%
\begin{array}{c}
  \Psi_1 \\
  \Psi_2 \\
\end{array}%
\right).$$

Here $u_1=\overline{u}_2,$ and the two scalar equations comprised in
(\ref{proj-eq}) are complex conjugated one to the other. Hence it
suffices to study one of the them. We denote $u_1=u$ and the projected equation  is
\begin{equation}\label{proj-eq1}
\frac{du}{dt}=(\mu+i\omega
)u+\Psi_1(0)\widetilde{f}(\varphi_1u+\varphi_2\overline{u}+\mathbf{v}).
\end{equation}

If the initial condition $\phi$ is taken on the center manifold,
then its image through the semigroup $\{T(t)\}_{t\geq 0},$
$T(\phi)=x_t(\phi)$, is still on the manifold. Hence
\begin{equation}\label{zt}T(t)(\phi)=u(t)\varphi_1+\overline{u(t)}\varphi_2+
w(u(t)\varphi_1+\overline{u(t)}\varphi_2),
\end{equation}
with $u(\cdot),$ solution of the equation

\begin{equation}\label{proj-eq1}
\frac{du}{dt}=\lambda_{1}u+\Psi_1(0)\widetilde{f}(\varphi_1u+\varphi_2\overline{u}+w(u\varphi_1+
\overline{u}\varphi_2)),
\end{equation}
with the initial condition $u(0)=u_0,$ where
$\mathcal{P}(\phi)=u_0\varphi_1+\overline{u}_0\varphi_2.$ The real
and the imaginary parts of this complex equation, represent the
two-dimensional restricted to the center manifold problem.

Now, if we use the series of powers of $w$ from \eqref{cm}, we can write
\[\widetilde{f}(\varphi_1u+\varphi_2\overline{u}+w(u\varphi_1+
\overline{u}\varphi_2))=\sum_{i+j\geq 2}\frac{1}{i!j!}f_{i,j}u^i\overline{u}^j,
\]
and, by putting
\begin{equation}\label{g-f}g_{i,j}=\Psi_1(0)f_{i,j},\end{equation}
eq. \eqref{proj-eq1} becomes
\begin{equation}\label{proj-eq2}
\frac{du}{dt}=\lambda_{1}u+\sum_{i+j\geq 2}\frac{1}{i!j!}g_{i,j}u^i\overline{u}^j.
\end{equation}
Remark that $f_{j,i}=\overline{f}_{i,j},$ since the function $\widetilde{f}$ has real values.

\section{Computation of the coefficients $w_{i,j}$}

For the determination of the functions $w_{i,j}$ of \eqref{cm}, the
following relation (obtained from the invariance of the center
manifold) is used \cite{WWPOG}, \cite{MNO}

\begin{equation}\label{gendiffeq}\frac{\partial}{\partial s}\sum_{j+k\geq2}\frac{1}{j!k!}w_{j,k}(s)u^{j}\overline{u}^{k}=
\sum_{j+k\geq2}\frac{1}{j!k!}g_{j,k}u^{j}\overline{u}^{k}\varphi_1(s)+
\end{equation}
\[+
\sum_{j+k\geq2}\frac{1}{j!k!}\overline{g}_{j,k}\overline{u}^{j}u^{k}\varphi_2(s)+\frac{\partial}{\partial
t}\sum_{j+k\geq2}\frac{1}{j!k!}w_{j,k}(s)u^{j}\overline{u}^{k}.
\]
This relation yields, by equating the terms of the same degree, differential equations for each $w_{i,j}.$ The integration constants are obtained from the following relation, also by equating the terms of same degree:

\begin{equation}\label{conddiffeq}
\frac{d}{dt}\sum_{j+k\geq2}\frac{1}{j!k!}w_{j,k}(0)u^{j}\overline{u}^{k}+
\sum_{j+k\geq2}\frac{1}{j!k!}g_{j,k}u^{j}\overline{u}^{k}\varphi_1(0)+
\sum_{j+k\geq2}\frac{1}{j!k!}\overline{g}_{j,k}\overline{u}^{j}u^{k}\varphi_2(0)=
\end{equation}
\[=A\sum_{j+k\geq2}\frac{1}{j!k!}w_{j,k}(0)u^{j}\overline{u}^{k}+B\sum_{j+k\geq2}\frac{1}{j!k!}w_{j,k}(-r)u^{j}\overline{u}^{k}+
\sum_{j+k\geq2}\frac{1}{j!k!}f_{j,k}u^{j}\overline{u}^{k}.
\]

\subsection{The system for  $w_{2,1}(-r),\,\,w_{2,1}(0)$}

By matching the terms that contain $u^2\overline{u},$  we get the
differential  equation and the condition for $w_{2,1}.$ These are

\[\frac{d }{ds}w_{2,1}(s)=\omega i w_{2,1}(s)+g_{2,1}e^{\omega i
s}+\overline{g}_{1,2}e^{-\omega is}+2w_{2,0}(s)g_{1,1}+\]
\[\quad\quad\quad\quad+w_{1,1}(s)g_{2,0}+2w_{1,1}(s)\overline{g}_{1,1}
+w_{0,2}(s)\overline{g}_{0,2},
\]
and
\[\omega iw_{2,1}(0)
+2w_{2,0}(0)g_{1,1}+w_{1,1}(0)g_{2,0}+2w_{1,1}(0)\overline{g}_{1,1}
+w_{0,2}(0)\overline{g}_{0,2}+g_{2,1}+\overline{g}_{1,2}=\]
\[
=Aw_{2,1}(0)+ B w_{2,1}(-r)+f_{2,1}.
\]
From these  we obtain the system of equations for $w_{2,1}(0)$  and
$w_{2,1}(-r)$:

\begin{equation}\label{alg-sys1}-e^{-\omega ir}w_{2,1}(0)+w_{2,1}(-r)=-g_{2,1}re^{-\omega
ir}+\frac{i}{2\omega} \overline{g}_{1,2}(e^{\omega ir}-e^{-\omega
ir})-\end{equation} \[-2g_{1,1}e^{-\omega
ir}\int_{-r}^0w_{2,0}(\theta)e^{-\omega i
\theta}d\theta-(g_{2,0}+2\overline{g}_{1,1})e^{-\omega
ir}\int_{-r}^0 w_{1,1}(\theta)e^{-\omega i\theta}d\theta -\]
\[-\overline{g}_{0,2}e^{-\omega ir}\int_{-r}^0
w_{0,2}(\theta)e^{-\omega i\theta}d\theta,
\]
\medskip

\[-(\omega i-A)w_{2,1}(0)+
Bw_{2,1}(-r)=g_{2,1}+\overline{g}_{1,2}-f_{2,1}+2g_{1,1}w_{2,0}(0)+
\]
\begin{equation}\label{alg-sys2}\quad\quad\quad\quad+
(g_{2,0}+2\overline{g}_{1,1})w_{1,1}(0)
+\overline{g}_{0,2}w_{0,2}(0).\end{equation}
The matrix of this
system is
\[
    \left(
\begin{array}{cc}
 -e^{-\omega ir}& 1\\
       -(\omega i-A)   & B \end{array}
     \right)
\]
and its determinant is  $\Delta:=\omega i-A-Be^{-\omega ir}$ that is
equal to zero, since $\omega i$ is a solution of the characteristic
equation \eqref{ch-eq}. Hence we obtained

\medskip

\textbf{Proposition 3.1.} \textit{The matrix of the system of
algebraic linear equations for} $w_{2,1}(0)$ \textit{and}
$w_{2,1}(-r)$ \textit{has null determinant}.

\medskip

In this situation the system has either an infinity of solutions or
no solutions at all. The following Proposition solves this problem.

\medskip

\textbf{Proposition 3.2.} \textit{The equations of the system
\eqref{alg-sys1}-\eqref{alg-sys2} are dependent.}

\smallskip

\textbf{Proof.} We denote the right hand side of \eqref{alg-sys1} by
$R_1$ and the right hand side of \eqref{alg-sys2} by $R_2.$ We have
to prove that
\begin{equation}\label{BR}BR_1=R_2.\end{equation} We
prove this by showing that the following relations hold:
$$B\left[-g_{2,1}re^{-\omega ir}+\frac{i}{2\omega}
\overline{g}_{1,2}(e^{\omega ir}-e^{-\omega
ir})\right]=g_{2,1}+\overline{g}_{1,2}-f_{2,1}; \eqno{(R1)} $$
$$-2Bg_{1,1}e^{-\omega
ir}\int_{-r}^0w_{2,0}(\theta)e^{-\omega
i\theta}d\theta=2g_{1,1}w_{2,0}(0); \eqno{(R2)}$$
$$-B(g_{2,0}+2\overline{g}_{1,1})e^{-\omega ir}\int_{-r}^0 w_{1,1}(\theta)e^{-\omega
i\theta}d\theta=(g_{2,0}+2\overline{g}_{1,1})w_{1,1}(0);
\eqno{(R3)}$$
$$-B\overline{g}_{0,2}e^{-\omega ir}\int_{-r}^0
w_{0,2}(\theta)e^{-\omega
i\theta}d\theta=\overline{g}_{0,2}w_{0,2}(0). \eqno{(R4)}$$

\vspace{0.5cm}

\textbf{(\textit{R}1).} We first replace $g_{2,1}$ with $\Psi_1(0)
f_{2,1},$ and $\overline{g}_{1,2}$ with $\overline{\Psi}_1(0)
f_{2,1}$  and (\textit{R}1) takes the form
\[B\left[-\Psi_{1}(0)f_{2,1}re^{-\omega ir}+\frac{i}{2\omega}
\overline{\Psi}_{1}(0)f_{2,1}(e^{\omega ir}-e^{-\omega
ir})\right]=\Psi_{1}(0)f_{2,1}+\overline{\Psi}_{1}(0)f_{2,1}-f_{2,1}.\]

 If $f_{2,1}=0,$ then  (\textit{R}1) is proved. If not,
after dividing the relation by $f_{2,1},$ we obtain
\begin{equation}\label{rel-Psi}-Be^{-\omega
ir}\Psi_1(0)r+B\frac{\overline{\Psi}_1(0)i}{2\omega }(e^{\omega
ir}-e^{-\omega ir})=\Psi_1(0) +\overline{\Psi}_1(0)-1.
\end{equation}

We prove that
\begin{equation}\label{rel-desc}-Be^{-\omega
ir}\Psi_1(0)r=\Psi_1(0)
-1,\;\mathrm{and}\;B\frac{\overline{\Psi}_1(0)i}{2\omega }(e^{\omega
ir}-e^{-\omega ir})=\overline{\Psi}_1(0).
\end{equation}

We write \eqref{rel-desc}$_1$ as
\[\Psi_1(0)+Be^{-\omega
ir}\Psi_1(0)r= 1
\]
and remark that it can be re-formulated as
\[\Psi_1(0)+B\int_{-r}^0\Psi_1(0)e^{-\omega
i(\theta+r)}e^{\omega i\theta}d\theta =1.
\]
But this last relation is nothing else than
\[\langle\Psi_1, \varphi_1 \rangle=1,
\]
that is true. Similarly, for \eqref{rel-desc}$_2$, we write it as
\[\overline{\Psi}_1(0)+B\frac{\overline{\Psi}_1(0)}{2\omega i
}(e^{\omega ir}-e^{-\omega ir})=0,
\]
and rewrite this as
\[\overline{\Psi}_1(0)+B\int_{-r}^0\overline{\Psi}_1(0)e^{\omega
i(\theta+r)}e^{\omega i\theta}d\theta=0,
\]
that is, actually,
\[\langle\Psi_2, \varphi_1 \rangle=0,
\]
that is true. Since the two equalities of \eqref{rel-desc} are
proved, by adding them, we obtain \eqref{rel-Psi}.

\vspace{0.5cm}

\textbf{Relations (\textit{R}2)-(\textit{R}4).} If $g_{1,1} =0,$
then relation (\textit{R}2) holds. We assume $g_{1,1} \neq 0.$
Similarly, for the two following relations, we assume
$g_{2,0}+2\overline{g}_{1,1}\neq 0$ and $\overline{g}_{0,2}\neq 0$.

We see that, after dividing with the assumed non-zero coefficients,
each of these relations has the form:
\begin{equation}\label{rel-j-k}-Be^{-\omega ir}\int_{-r}^0w_{j,k}(\theta)e^{-\omega
i\theta}d\theta=w_{j,k}(0)
\end{equation}
($j,k\geq 0,\,j+k=2$) that can be written as
\begin{equation}\label{rel-j-k1}w_{j,k}(0)+B\int_{-r}^0e^{-\omega
i(\theta+r)}w_{j,k}(\theta)d\theta=0.
\end{equation}
But relation \eqref{rel-j-k1} is in fact
\[\langle \psi_2, w_{j,k}\rangle =0,
\]
which is true, since each of the functions $w_{j,k}$ belong to the
complementary of $\mathcal{M}.$ Hence all relations (\textit{R}2) -
(\textit{R}4) are proved.

\smallskip

 Since we proved relations (\textit{R}1) - (\textit{R}4),
by adding them, relation \eqref{BR} follows, and with it, the
conclusion of our Proposition. $\Box$

The problem that occurs at this point is how to choose some "proper"
values of $w_{2,1}(0),\,w_{2,1}(-r)$ from the infinity of solutions
of the system \eqref{alg-sys1} - \eqref{alg-sys2}. In the next
section we give a solution.

\section{How to compute $w_{2,1}$}

We  consider a perturbation of our problem, of the form

\begin{equation}\label{pert-eq}\dot{x}(t)=A_\epsilon x(t)+B_\epsilon x(t-r)+\widehat{f}(x(t),\,x(t-r)),
\end{equation}
where $A_\epsilon, \,B_\epsilon$ depend smoothly enough on $\epsilon>0,$
$\displaystyle\lim_{\epsilon\searrow 0}A_\epsilon = A,\, \lim_{\epsilon\searrow 0}B_\epsilon = B,$
and are chosen such that, for small enough $\epsilon$,  the linearized problem attached to \eqref{pert-eq} admits the eigenvalues $\lambda_{\epsilon 1,2}=\mu_\epsilon\pm \omega_\epsilon i,$ with $\mu_\epsilon>0,$ while all other eigenvalues have negative real part.

Obviously, from the construction it follows that
$\displaystyle\lim_{\epsilon\searrow 0}\mu_\epsilon = 0,\,
\lim_{\epsilon\searrow 0}\omega_\epsilon = \omega$. An example of
perturbed problem is given at the end of this subsection.

From the construction, it follows that problem \eqref{pert-eq}
admits an unstable manifold, tangent to the space
$\mathcal{M}_\epsilon$ spanned by the two eigenfunctions
$\varphi_{\epsilon 1,2}(s)=e^{(\mu_\epsilon\pm i\omega_\epsilon)s},$
corresponding to the two eigenvalues $\lambda_{\epsilon 1,2}.$

The unstable manifold is the graph of a function $w_{\epsilon}$ defined on $\mathcal{M}_\epsilon$ and taking values in a subspace $\mathcal{N}_\epsilon$ of $\mathcal{B}_C,$ complementary to $\mathcal{M}_\epsilon$.

A procedure, similar to that presented in Section 2 for the non-perturbed problem, is developed in order to construct the restriction of the problem to the unstable manifold:

-the adjoint equation and the eigenfunctions of its linearized, i.e. the functions $\psi_{\epsilon 1}(s)=e^{-\lambda_\epsilon s},\,\psi_{\epsilon 2}(s)=e^{-\overline{\lambda}_\epsilon s},\,s\in [-r,0],$ are considered;

-the corresponding bilinear form, denoted  also by $\langle\, \cdot\, ,\,\cdot\,\rangle$ is constructed,  and the functions $\Psi_{\epsilon 1},\,\Psi_{\epsilon 2}\in C([0,r],\mathbb{C})$ are computed such that $\langle \Psi_{\epsilon j},\varphi_{\epsilon k} \rangle =\delta_{j,k};$ these are
\[\Psi_{\epsilon 1}(s)=\frac{
1+(\overline{\lambda}_\epsilon- A_{\epsilon})r}{(1-A_{\epsilon}r+\mu_\epsilon
r)^2+\omega_\epsilon^2r^2}\psi_{\epsilon 1},\;\Psi_{\epsilon 2}=\overline{\Psi}_{\epsilon 1};
\]

- the projector $\mathcal{P}_\epsilon:\mathcal{B}_C\rightarrow \mathcal{M}_\epsilon,$ such that $I-\mathcal{P}_\epsilon:\mathcal{B}_C\rightarrow \mathcal{N}_\epsilon$, is defined by
$$\mathcal{P}_\epsilon(\phi)=\langle \Psi_{\epsilon 1},\phi\rangle \varphi_1+\langle \Psi_{\epsilon 2},\phi\rangle \varphi_2.$$

We remark that, by the construction of the bilinear form, and by the definition of the function $w_\epsilon$ whose graph is the unstable manifold, we have:
\[\langle \Psi_{\epsilon j}, w_\epsilon \rangle=0,\,\,j=1,2.
\]

With the same type of reasonings (coming from \cite{Far}), we find
that the problem reduced to the unstable manifold is
\[\frac{dv}{dt}=\lambda_{\epsilon}v+\Psi_{\epsilon 1}(0)\widetilde{f}(v\varphi_{\epsilon 1}+\overline{v}\varphi_{\epsilon 2}+w_\epsilon(v\varphi_{\epsilon 1}+\overline{v}\varphi_{\epsilon
2})).
\]

As we did in Section 1 for $w,$ we consider the function
$$\widehat{w}_\epsilon(v,
\overline{v}):=w_\epsilon(v\varphi_{\epsilon
1}+\overline{v}\varphi_{\epsilon 2})$$ and (dropping the hat for
simplicity of notations), we write:
\begin{equation}\label{um}  w_\epsilon(v,\overline{v})=\sum_{i+j\geq 2}\frac{1}{i!j!}w_{\epsilon i,j}v^i\overline{v}^j.\end{equation}
The coefficients  $w_{\epsilon i,j}$ are found by solving
differential equations coming from  relations similar to
\eqref{gendiffeq}, \eqref{conddiffeq} with $w_{j,k}, \,g_{j,k},...$
replaced by  $w_{\epsilon j,k}, \,g_{\epsilon j,k},...$.

The equation for $w_{\epsilon 2,1}$ is
\[\frac{d w_{\epsilon2,1}(s)}{ds}=(2\lambda_\epsilon+\overline{\lambda}_\epsilon) w_{\epsilon 2,1}(s)+g_{\epsilon 2,1}e^{\lambda_\epsilon
s}+\overline{g}_{\epsilon 1,2}e^{\overline{\lambda}_\epsilon s}
+2w_{\epsilon 2,0}(s)g_{\epsilon 1,1}+\]
\[+w_{\epsilon 1,1}(s)g_{\epsilon 2,0}+2w_{\epsilon 1,1}(s)\overline{g}_{\epsilon 1,1}
+w_{\epsilon 0,2}(s)\overline{g}_{\epsilon 0,2},
\]
while the condition to determine the integration constant is
\[(2\lambda_\epsilon+\overline{\lambda}_\epsilon) w_{\epsilon 2,1}(0)
+2w_{\epsilon 2,0}(0)g_{11}+w_{\epsilon
1,1}(0)g_{\epsilon2,0}+2w_{\epsilon1,1}(0)\overline{g}_{\epsilon1,1}
+w_{\epsilon0,2}(0)\overline{g}_{\epsilon0,2}+g_{\epsilon2,1}+
\overline{g}_{\epsilon1,2}=\]
\[
=A_{\epsilon}w_{\epsilon2,1}(0)+
B_{\epsilon}w_{\epsilon2,1}(-r)+f_{\epsilon2,1}.
\]
By integrating the differential equation above between $-r$ and $0$, we find the system of equations for
$w_{\epsilon2,1}(-r),\,w_{\epsilon2,1}(0):$
\begin{equation}\label{sist-eps1}-e^{-(2\lambda_\epsilon+\overline{\lambda}_\epsilon)r}w_{\epsilon2,1}(0)+w_{\epsilon2,1}(-r)=
\frac{-1}{\lambda_\epsilon+\overline{\lambda}_\epsilon}g_{\epsilon2,1}(
e^{-\lambda_\epsilon
r}-e^{-(2\lambda_\epsilon+\overline{\lambda}_\epsilon)r})-\;\;\;
\end{equation}
\[-\frac{1}{2\lambda_\epsilon}\overline{g}_{\epsilon1,2}(e^{-\overline{\lambda}_\epsilon
r}-e^{-(2\lambda_\epsilon+\overline{\lambda}_\epsilon)r})-2g_{\epsilon1,1}e^{-(2\lambda_\epsilon+\overline{\lambda}_\epsilon)r}
\int_{-r}^0w_{\epsilon2,0}(\theta)e^{-(2\lambda_\epsilon+
\overline{\lambda_\epsilon}) \theta}d\theta-\]

\[-(g_{\epsilon2,0}+2\overline{g}_{\epsilon1,1})e^{-(2\lambda_\epsilon+\overline{\lambda}_\epsilon)r}
\int_{-r}^0
w_{\epsilon1,1}(\theta)e^{-(2\lambda_\epsilon+\overline{\lambda}_\epsilon)
\theta}d\theta-
\]
\[-\overline{g}_{\epsilon0,2}e^{-(2\lambda_\epsilon+\overline{\lambda}_\epsilon)r}\int_{-r}^0
w_{\epsilon0,2}(\theta)
e^{-(2\lambda_\varepsilon+\overline{\lambda}_\epsilon)
\theta}d\theta,
\]
\vspace{0.8cm}
\begin{equation}\label{sist-eps2}
\left(A_{\epsilon}-2\lambda_\epsilon-\overline{\lambda}_\epsilon\right)w_{\epsilon2,1}(0)+
B_{\epsilon}w_{\epsilon2,1}(-r)=g_{\epsilon2,1}+
\overline{g}_{\epsilon1,2}-f_{\epsilon2,1}+
\end{equation}
\[+2w_{\epsilon 2,0}(0)g_{11}+w_{\epsilon
1,1}(0)g_{\epsilon2,0}+2w_{\epsilon1,1}(0)\overline{g}_{\epsilon1,1}
+w_{\epsilon0,2}(0)\overline{g}_{\epsilon0,2}.
\]

The matrix of the system is
\[
\left(%
\begin{array}{cc}
  -e^{-(2\lambda_\epsilon+\overline{\lambda}_\epsilon)r} & 1 \\
  A_{\epsilon}-2\lambda_\epsilon-\overline{\lambda}_\epsilon & B_{\epsilon} \\
\end{array}%
\right)
\]
with determinant
\begin{equation}\label{Delta-eps}\Delta_\epsilon=-B_{\epsilon}e^{-(2\lambda_\epsilon+\overline{\lambda}_\epsilon)r}-
A_{\epsilon}+2\lambda_\epsilon+\overline{\lambda}_\epsilon,\end{equation}
that is different of zero because otherwise the number
$2\lambda_\epsilon+\overline{\lambda}_\epsilon$ would be an
eigenvalue, with real part equal to $3\mu_\epsilon$, that
contradicts the fact that all eigenvalues have real part $\leq
\mu_\epsilon$.

We denote the right-hand side of the
two equations by $R_{\epsilon1},\,R_{\epsilon2}, $ respectively.

\medskip

\textbf{Proposition 4.1} \emph{When $\epsilon\rightarrow 0,$ the
coefficients of system \eqref{sist-eps1}-\eqref{sist-eps2} tend to
the coefficients of system \eqref{alg-sys1}-\eqref{alg-sys2}.}

\smallskip

\textbf{Proof.} The assertion of the proposition is obvious for all
the coefficients of the unknowns and for the terms in the right-hand
sides, excepting the first term from the right hand sides of
\eqref{alg-sys1} and \eqref{sist-eps1}.

But also for these the conclusion comes easily
\[\lim_{\epsilon\rightarrow 0} \frac{-1}{\lambda_\epsilon+\overline{\lambda}_\epsilon}
g_{\epsilon2,1}( e^{-\lambda_\epsilon
r}-e^{-(2\lambda_\epsilon+\overline{\lambda}_\epsilon)r})=\]
\[=-\lim_{\epsilon\rightarrow
0} g_{\epsilon2,1}e^{-\lambda_\epsilon r}\lim_{\epsilon\rightarrow
0}\frac{1-e^{-(\lambda_\epsilon+\overline{\lambda}_\epsilon)r}}{(\lambda_\epsilon+\overline{\lambda}_\epsilon)r}r
 =-g_{2,1}e^{-\omega i r}r,
\]
since $\lambda_\epsilon+\overline{\lambda}_\epsilon=2\mu(\epsilon)\rightarrow 0,$
when $\epsilon\rightarrow 0. \,\,\Box$

\vspace{0.5cm}

Now, the natural idea is to solve system
\eqref{sist-eps1}-\eqref{sist-eps2} and to compute the limit of its
solution when $\epsilon\rightarrow 0.$ As a matter of fact, it is
enough to find $w_{\epsilon2,1}(0),$ to compute its limit, and,
then, $w_{2,1}(-r)$ will be found from one of the equations
\eqref{alg-sys1}, \eqref{alg-sys2}. By denoting the right hand sides
of eqs. \eqref{sist-eps1}, \eqref{sist-eps2} with
$R_{\epsilon1},\,R_{\epsilon2},$ respectively, we have

\[
    w_{\epsilon2,1}(0) =\frac{B_{\epsilon}R_{\epsilon1}-R_{\epsilon2}}{\Delta_\epsilon}.
    \]

 We know that both the numerator and the denominator
of the above expressions tend to 0 when $\epsilon\rightarrow 0$, but
the next Proposition shows that we can overcome this problem.

\textbf{Proposition 4.2} \textit{For any $\epsilon>0$,
$B_{\epsilon}R_{\epsilon1}-R_{\epsilon2}$ and }$\Delta_\epsilon$
\textit{can be written in the form }
\[B_{\epsilon}R_{\epsilon1}-R_{\epsilon2}=\mu_\epsilon
h_1(\epsilon),
\]
\[\Delta_\epsilon=\mu_\epsilon
h_2(\epsilon),
\]
\textit{where the functions $h_j,\,j=1,2$ have finite limit for
$\epsilon\rightarrow 0,$ and $\displaystyle
\lim_{\epsilon\rightarrow 0}h_2(\epsilon)\neq 0.$  }

\textbf{Proof.} First we deal with the determinant
$\Delta_\epsilon.$ We have, by using the characteristic equation
($\lambda_\epsilon=A_\epsilon+B_\epsilon e^{-\lambda_\epsilon r}$),
\[\Delta_\epsilon=-B_{\epsilon}e^{-(2\lambda_\epsilon+\overline{\lambda}_\epsilon)r}+B_\epsilon e^{-\lambda_\epsilon r}+\lambda_\epsilon+\overline{\lambda}_\epsilon=
\]
\[=B_\epsilon e^{-\lambda_\epsilon r}\left(1-e^{-2\mu_\epsilon r}\right)+2\mu_\epsilon=
B_\epsilon e^{-\lambda_\epsilon r}2\mu_\epsilon r\left(1-\frac{2\mu_\epsilon r}{2!}+\frac{(2\mu_\epsilon r)^2}{3!}+...\right)+2\mu_\epsilon=
\]
\[=2\mu_\epsilon\left[B_\epsilon e^{-\lambda_\epsilon r} r\left(1-\frac{2\mu_\epsilon r}{2!}+\frac{(2\mu_\epsilon r)^2}{3!}+...\right)+1\right].
\]
Hence the assertion concerning $\Delta_\epsilon$ holds, the
expression of $h_2$ being obvious from the above relations.

Moreover, we see that $\displaystyle \lim_{\epsilon\rightarrow
0}h_2(\epsilon)=2\left(Be^{-\omega ir }r+1\right)=2r\omega i
-2rA+2,$ that can not be zero, since $\omega \neq 0.$

Now, in order to treat the term $B_{\epsilon}R_{\epsilon1}-R_{\epsilon2}$, inspired by the proof of Proposition 3.2, we write
\[B_{\epsilon}R_{\epsilon1}-R_{\epsilon2}=E_1+E_2+E_3+E_4,
\]
where
\[E_1=\frac{-B_\epsilon}{\lambda_\epsilon+\overline{\lambda}_\epsilon}g_{\epsilon2,1}(
e^{-\lambda_\epsilon
r}-e^{-(2\lambda_\epsilon+\overline{\lambda}_\epsilon)r})+
\frac{-B_\epsilon}{2\lambda_\epsilon}\overline{g}_{\epsilon1,2}(e^{-\overline{\lambda}_\epsilon
r}-e^{-(2\lambda_\epsilon+\overline{\lambda}_\epsilon)r})-(g_{\epsilon2,1}+
\overline{g}_{\epsilon1,2}-f_{\epsilon2,1}),
\]
\[E_2=-2B_\epsilon g_{\epsilon1,1}e^{-(2\lambda_\epsilon+\overline{\lambda}_\epsilon)r}
\int_{-r}^0w_{\epsilon2,0}(\theta)e^{-(2\lambda_\epsilon+
\overline{\lambda_\epsilon}) \theta}d\theta-2g_{\epsilon1,1}w_{\epsilon 2,0}(0),
\]
\[E_3=-B_\epsilon(g_{\epsilon2,0}+2\overline{g}_{\epsilon1,1})e^{-(2\lambda_\epsilon+\overline{\lambda}_\epsilon)r}
\int_{-r}^0
w_{\epsilon1,1}(\theta)e^{-(2\lambda_\epsilon+\overline{\lambda}_\epsilon)
\theta}d\theta-(g_{\epsilon2,0}+2\overline{g}_{\epsilon1,1})w_{\epsilon
1,1}(0),
\]
\[E_4=-B_\epsilon\overline{g}_{\epsilon0,2}e^{-(2\lambda_\epsilon+\overline{\lambda}_\epsilon)r}\int_{-r}^0
w_{\epsilon0,2}(\theta)
e^{-(2\lambda_\varepsilon+\overline{\lambda}_\epsilon)
\theta}d\theta-\overline{g}_{\epsilon0,2}w_{\epsilon0,2}(0).
\]
$\mathbf{E_1}.$ We write  $E_1$ as
$E_1=f_{\epsilon2,1}(E_{11}+E_{12}),$ where:
\[E_{11}=\frac{-B_\epsilon}{\lambda_\epsilon+\overline{\lambda}_\epsilon}\Psi_{\epsilon1}(0)(
e^{-\lambda_\epsilon
r}-e^{-(2\lambda_\epsilon+\overline{\lambda}_\epsilon)r})
-\Psi_{\epsilon1}(0)
+1,
\]
\[E_{12}=\frac{-B_\epsilon}{2\lambda_\epsilon}\overline{\Psi}_{\epsilon1}(0)(e^{-\overline{\lambda}_\epsilon
r}-e^{-(2\lambda_\epsilon+\overline{\lambda}_\epsilon)r})-\overline{\Psi}_{\epsilon1}(0).
\]
For $E_{11},$  we have
\[E_{11}=1-\left(\Psi_{\epsilon1}(0)+B_\epsilon \Psi_{\epsilon1}(0)e^{-\lambda_\epsilon r}\int_{-r}^0e^{(\lambda_\epsilon+\overline{\lambda}_\epsilon)s}ds\right)=\]
\[=1-\left(\Psi_{\epsilon1}(0)+B_\epsilon \Psi_{\epsilon1}(0)\int_{-r}^0e^{-\lambda_\epsilon(s+ r)}e^{(2\lambda_\epsilon+\overline{\lambda}_\epsilon)s}ds\right)=
\]
\[=1-\langle \Psi_{\epsilon1}, \eta_\epsilon \rangle,
\]
where
$\eta_\epsilon(s)=e^{(2\lambda_\epsilon+\overline{\lambda}_\epsilon)s},\,\,s\in[-r,0].$

But we know that $1=\langle  \Psi_{\epsilon1}, \varphi_{\epsilon 1}\rangle,$ and, substituting this above, we find
\[E_{11}=\langle  \Psi_{\epsilon1}, \varphi_{\epsilon1}\rangle-\langle \Psi_{\epsilon1}, \eta_\epsilon \rangle=
\langle  \Psi_{\epsilon1}, \varphi_{\epsilon1}-\eta_\epsilon\rangle.
\]
But
\begin{equation}\label{imp-rel}\varphi_{\epsilon1}(s)-\eta_\epsilon(s)=e^{\lambda_\epsilon s}-e^{(2\lambda_\epsilon+\overline{\lambda}_\epsilon)s}=e^{\lambda_\epsilon s}\left(1- e^{2\mu_\epsilon s}\right)= \end{equation}
\[=-2\mu_\epsilon s\, e^{\lambda_\epsilon s}\left(1+\frac{2\mu_\epsilon s}{2!}+\frac{(2\mu_\epsilon s)^2}{3!}+...\right).
\]
The series in the paranthesis is convergent and its sum is a bounded
function on $[-r,0].$  We set
\[\rho_\epsilon(s)=-2s\, e^{\lambda s}\left(1+\frac{2\mu_\epsilon s}{2!}+\frac{(2\mu_\epsilon s)^2}{3!}+...\right)
\]
and we obtain
\[E_{11}=\mu_\epsilon \langle \Psi_{\epsilon1}, \,\rho_\epsilon  \rangle.
\]

Now we pass to $E_{12}$
\[E_{12}=-\overline{\Psi}_{\epsilon1}(0)-\frac{B_\epsilon}{2\lambda_\epsilon}\overline{\Psi}_{\epsilon1}(0)(e^{-\overline{\lambda}_\epsilon
r}-e^{-(2\lambda_\epsilon+\overline{\lambda}_\epsilon)r})=
\]
\[=-\overline{\Psi}_{\epsilon1}(0)-B_\epsilon\overline{\Psi}_{\epsilon1}(0)e^{-\overline{\lambda}_\epsilon
r}\int_{-r}^0e^{2\lambda_\epsilon s}ds=
\]
\[=-\overline{\Psi}_{\epsilon1}(0)-B_\epsilon\overline{\Psi}_{\epsilon1}(0)\int_{-r}^0e^{-\overline{\lambda}_\epsilon (s+r)}e^{(2\lambda_\epsilon+\overline{\lambda}_\epsilon )s}ds=
-\langle \Psi_{\epsilon2}, \eta_\epsilon  \rangle.
\]
But  $\langle \Psi_{\epsilon2}, \varphi_{\epsilon1}\rangle=0$, and, then, we can write
\[E_{12}= \langle \Psi_{\epsilon2}, \varphi_{\epsilon1}\rangle-\langle \Psi_{\epsilon2}, \eta_\epsilon  \rangle= \langle \Psi_{\epsilon2}, \varphi_{\epsilon1}-\eta_\epsilon\rangle,
\]
and from this point, by repeating step by step the reasonings made for $E_{11}$ we obtain
\[E_{12}=\mu_\epsilon \langle \Psi_{\epsilon2}, \,\rho_\epsilon  \rangle.
\]
Finally,
\begin{equation}E_{1}=\mu_\epsilon f_{\epsilon2,1}\left( \langle \Psi_{\epsilon1}, \,\rho_\epsilon  \rangle+\langle \Psi_{\epsilon2}, \,\rho_\epsilon  \rangle\right).
\end{equation}

\noindent$\mathbf{E_2}\,-\,\mathbf{E_4}.$ For a unitary writing, we
define
\[\alpha_{2,0}:=-2g_{\epsilon1,1},\,\,\alpha_{1,1}:=-g_{\epsilon2,0}-2\overline{g}_{\epsilon1,1},\,
\,\alpha_{0,2}:=-\overline{g}_{\epsilon0,2}.\]

Then, each of the expressions $E_i,\, i=2,3,4, $ can be written as
\[\alpha_{j,k}\left(w_{\epsilon j,k}(0)+ B_\epsilon \int_{-r}^0 e^{-(2\lambda_\epsilon+\overline{\lambda}_\epsilon)(s+r)} w_{\epsilon j,k}(s)ds     \right)=\alpha_{j,k}\langle \widetilde{\eta}_\epsilon, w_{\epsilon j,k}\rangle,\;\;
\]
where $j,k>0,\,j+k=2,$ and
$\widetilde{\eta}_\epsilon(\zeta)=e^{-(2\lambda_\epsilon+\overline{\lambda}_\epsilon)\zeta},\,\,
\zeta\in [0,r].$

From the definitions of the function $w_\epsilon,$ that defines the
invariant manifold, and of the projector on $\mathcal{M}_\epsilon,$
we have
\[\langle \psi_{\epsilon1}, w_{\epsilon j,k} \rangle =0.
\]
Then we may write
\[\langle \widetilde{\eta}_\epsilon, w_{\epsilon j,k}\rangle=
\langle \widetilde{\eta}_\epsilon, w_{\epsilon j,k}\rangle-\langle
\psi_{\epsilon1}, w_{\epsilon j,k} \rangle= \langle
\widetilde{\eta}_\epsilon-\psi_{\epsilon1}, w_{\epsilon j,k}\rangle,
\]
and
\[\widetilde{\eta}_\epsilon(\zeta)-\psi_{\epsilon1}(\zeta)=e^{-(2\lambda_\epsilon+\overline{\lambda}_\epsilon)\zeta}-e^{-\lambda_\epsilon\zeta}=e^{-\lambda_\epsilon\zeta}\left[ e^{-2\mu_\epsilon\zeta}-1 \right]=
\]
\[=-2\mu_\epsilon\zeta e^{-\lambda_\epsilon\zeta}\left( 1-\frac{2\mu_\epsilon\zeta}{2!}+\frac{(2\mu_\epsilon\zeta)^2}{3!}-... \right).
\]
We consider the function $\widetilde{\rho}:[0,r]\mapsto\mathbb{R},$
\[\widetilde{\rho}_\epsilon(\zeta)=-2\zeta e^{-\lambda_\epsilon\zeta}\left( 1-\frac{2\mu_\epsilon\zeta}{2!}+\frac{(2\mu_\epsilon\zeta)^2}{3!}-...
\right),
\]
and, finally we can write
$B_{\epsilon}R_{\epsilon1}-R_{\epsilon2}=\mu_\epsilon h_1(\epsilon),
$ where
\[\begin{array}{ccc}
    h_1(\epsilon) & = & f_{\epsilon2,1}\left( \langle \Psi_{\epsilon1},
\,\rho_\epsilon \rangle+\langle \Psi_{\epsilon2}, \,\rho_\epsilon
\rangle\right)-2g_{\epsilon1,1}\langle \widetilde{\rho}_\epsilon,
w_{\epsilon2,0}\rangle- \\
\,&\,&\,\\
    \, & \, &  -(g_{\epsilon2,0}+2\overline{g}_{\epsilon1,1})\langle
\widetilde{\rho}_\epsilon, w_{\epsilon1,1}\rangle
-\overline{g}_{\epsilon0,2}\langle \widetilde{\rho}_\epsilon,
w_{\epsilon0,2}\rangle.\\
  \end{array}
\]
The proof is complete. $\square$

\medskip

 Now we can easily compute
$w_{2,1}(0).$

\medskip

\noindent\textbf{Proposition 4.3.} \textit{The value of}
$w_{2,1}(0)$ \textit{is}
\begin{equation}\label{w21}w_{2,1}(0)=
\frac{f_{2,1}\langle \Psi_{1}+\Psi_{2}, \,\rho
\rangle-2g_{1,1}\langle \widetilde{\rho}, w_{2,0}\rangle
-(g_{2,0}+2\overline{g}_{1,1})\langle \widetilde{\rho},
w_{1,1}\rangle -\overline{g}_{0,2}\langle \widetilde{\rho},
w_{0,2}\rangle }{2r\omega i -2rA+2},
\end{equation}
\textit{where} $\rho(s)=-2se^{\omega is}, \;s\in [-r,0]$
\textit{and} $\widetilde{\rho}(\zeta)=-2\zeta e^{-\omega i\zeta},
\;\zeta\in [0,r].$
\medskip

\noindent\textbf{Proof.} The preceding Proposition implies
\[\lim_{\epsilon \rightarrow 0}w_{\epsilon2,1}(0)=
\lim_{\epsilon \rightarrow 0}\frac{h_1(\epsilon)}{h_2(\epsilon)}.
\]
In order to see that the coefficients $g_{\epsilon j,k}$ and the
functions $w_{\epsilon j,k}$ in $h_1(\epsilon)$ do not present any
problem when passing to limit, we list their values in the Appendix.

We see that $\lim_{\epsilon \rightarrow 0}g_{\epsilon j,k}=g_{j,k}$,
and that $w_{\epsilon j,k}\rightarrow w_{j,k}$ when
$\epsilon\rightarrow 0,$ uniformly on $[-r,0].$ Then, by observing
that
\[\lim_{\epsilon \rightarrow 0}\rho_\epsilon (s)=-2se^{\omega is},
\]
the convergence being uniform with respect to  $s\in[-r,0]$ and
\[\lim_{\epsilon \rightarrow 0}\widetilde{\rho}_\epsilon (\zeta)=-2\zeta e^{-\omega i\zeta},
\]
(uniform convergence on $[0,r]$), the result of our Proposition is
obtained. $\Box$

\bigskip

From any of the two equations of system \eqref{alg-sys1} -
\eqref{alg-sys2} we then find $w_{2,1}(-r).$

From the proof of Proposition 4.2, we see that $\displaystyle
\lim_{\epsilon\rightarrow 0}w_{\epsilon2,1}(0)$ does not depend on
the specific perturbation of the problem chosen.

\medskip
\textbf{Example of a concrete perturbation.} We take
$B_\epsilon=B(1+\epsilon)$ and chose $A_\epsilon$ in order to have
two eigenvalues  of the form $\mu(\epsilon)\pm \omega i$, with
$\mu(\epsilon)>0.$

 The characteristic equation associated to the
linear part of the equation is
\begin{equation}\label{char-pert}\lambda=A_\epsilon+B(1+\epsilon)e^{-\lambda r}.
\end{equation}
 Taking the imaginary parts of the equation, we find
\begin{equation}\label{miu-eps}\mu(\epsilon)=\ln\left[-\frac{B\sin(\omega r)}{\omega}(1+\epsilon)\right].
\end{equation}
Since $\omega i$ is an eigenvalue for the linearized of our problem,
the characteristic equation for $\epsilon=0$ implies
\[-\frac{B\sin(\omega r)}{\omega}=1,
\]
and, thus, for $\epsilon>0,$ $\displaystyle-\frac{B\sin(\omega
r)}{\omega}(1+\epsilon)>1$ and $\mu(\epsilon)>0.$

By taking the real part of \eqref{char-pert}, we obtain
\begin{equation}\label{A1}A_{\epsilon}=\mu(\epsilon)-B(1+\epsilon)e^{-\mu(\epsilon)}\cos(\omega r).
\end{equation}

Problem \eqref{pert-eq} with $A_\epsilon,\,B_\epsilon$ given above
is an example of perturbed problem. For small enough $\epsilon,$ its
linearized part has the eigenvalues $\lambda_{\epsilon 1,2
}=\mu(\epsilon)\pm i\omega,$ with $\mu(\epsilon)>0,$ while all other
eigenvalues have negative real part.

\medskip

\textbf{Remarks. 1.} In \cite{WWPOG}, for a particular problem, when
trying to  compute $w_{2,1}$ (denoted otherwise there),  a perturbed
problem (depending on a small parameter $\epsilon$) is considered,
the corresponding $w_{\epsilon2,1}$ is computed and the limit when
$\epsilon\rightarrow 0$ is taken. But, besides the fact that it is a
particular problem (hence the result is not general), there is no
proof there that the limit does not depend on the
specific perturbation chosen. \\

\textbf{2.} In a previous work \cite{AVI}, in which we intended to
present an example of Bautin type bifurcation in a delay
differential equation, we encountered the lack of uniqueness in
determining $w_{2,1}(0),\,w_{2,1}(-r).$ The equation considered
there, that did not come from some model (it was constructed by us)
was
\begin{equation}\label{ex}\dot{x}=ax(t-r)+x^2(t)+cx(t)x(t-r),
\end{equation}
with $r=\pi/2.$ Obviously, $x=0$ is an equilibrium point, and we
looked for the values of the parameters $a,\,c$ where the sufficient
conditions for Bautin type bifurcation are fulfilled \cite{AVI-0}.
We showed that the linear problem associated has, for $a_0=-1$ two
pure imaginary eigenvalues, $\lambda_{1,2}=\pm i,$ and that for two
values of $c,$
\[c_{1,2}=\frac{18-7\pi\pm\sqrt{36+212\pi+\pi^2}}{2(3\pi-2)},
\]
the first Lyapunov coefficient is zero. Up to this point the values
of $w_{2,1}$ were not necessary.

Then, we intended to compute the second Lyapunov coefficient since,
as we remarked in \cite {AVI-0}, when this is positive, the two
limit cycles (one inside the other) that occur in the Bautin type
bifurcation, exist when two eigenvalues with positive real part
exist, i.e. the cycles actually exist on the unstable manifold.

When computing the second Lyapunov coefficient $l_2(c)$, we needed
the values of $w_{2,1}(0),\,w_{2,1}(-r),$ but since the two
equations that yield these values are dependent, at that moment we
have chosen arbitrarily $w_{2,1}(0)=0$ and computed $w_{2,1}(-r)$ from one of
the two equations. We found that $l_2(c_1)>0,$ while $l_2(c_2)<0.$

However, the correct way of solving the problem is the one presented
in this paper, and this Remark is intended to play the role of an
Erratum to \cite{AVI}.

With the method developed  here, by using formula \eqref{w21} we find,  in the case of
$c_1(\approx1.52799),$ $w_{\epsilon2,1}(0)=0.285-0.28 i,\,w_{\epsilon2,1}(-r)=1.442-1.612 i,\,l_2=4.528.$
Hence, in the
parameters plane, in a neighborhood of the point $a=-1,\,c=c_1,$  there is a zone where the unstable manifold exists
and for parameters $a,c$ in a subset of this zone, two periodic
orbits (one inside the other) exist on the unstable manifold.

We re-analyzed the case of $c_2(\approx-2.06554)$ and found:\\
$w_{\epsilon2,1}(0)=-0.69-0.278 i,\,w_{\epsilon2,1}(-r)=-4.732-1.537 i,\,l_2=3.726.$

This shows that equation \eqref{ex} presents Bautin type bifurcation
for both values $c_{1},\, c_2$.

\section{Appendix}

First we write down the expressions of $f_{\epsilon i,j},\,i+j=2,$
and that of $f_{\epsilon 2,1}.$ We have:
\[f_{\epsilon2,0}=C_{2,0}+2C_{1,1}e^{-\lambda_\epsilon r}+C_{0,2}e^{-2\lambda_\epsilon r},
\]
\[f_{\epsilon1,1}=C_{2,0}+C_{1,1}(e^{-\lambda_\epsilon r}+e^{-\overline{\lambda}_\epsilon r})+
C_{0,2}e^{-2\mu_\epsilon r},
\]
\[f_{\epsilon0,2}=\overline{f_{\epsilon2,0}}=C_{2,0}+2C_{1,1}e^{-\overline{\lambda}_\epsilon r}+
C_{0,2}e^{-2\overline{\lambda}_\epsilon r},
\]
\[f_{\epsilon2,1}=C_{2,0}\left(2w_{\epsilon1,1}(0)+w_{\epsilon2,0}(0)
\right)+\]
\[+C_{1,1}\left(w_{\epsilon2,0}(0)e^{-\overline{\lambda}_\epsilon r}+
2w_{\epsilon1,1}(0)e^{-\lambda_\epsilon
r}+w_{\epsilon2,0}(-r)+2w_{\epsilon1,1}(-r)\right)+
\]
\[+C_{0,2}\left(2w_{\epsilon1,1}(-r)e^{-\lambda_\epsilon r}+w_{\epsilon2,0}(-r)
e^{-\overline{\lambda}_\epsilon r}\right)+C_{3,0}+
C_{2,1}\left(2e^{-\lambda_\epsilon
r}+e^{-\overline{\lambda}_\epsilon r} \right)+
\]
\[+C_{1,2}\left(2e^{-2\mu_\epsilon r}+e^{-2\lambda_\epsilon r}\right)+
C_{0,3}\left(e^{-\overline{\lambda}_\epsilon r}e^{-2\lambda_\epsilon
r}\right).
\]
By multiplying the above quantities by $\Psi_{\epsilon1}(0)$
 we obtain the corresponding $g_{\epsilon i,j}.$

Now we compute $w_{\epsilon j,k},\,j+k=2.$ The differential equation
for $w_{\epsilon2,0}$ is
\[
\frac{d w_{\epsilon2,0}(s)}{ds}=2\lambda_\epsilon
w_{\epsilon2,0}(s)+ g_{\epsilon2,0}e^{\lambda_\epsilon
s}+\overline{g}_{\epsilon0,2}e^{\overline{\lambda}_\epsilon s},
\]
and by integrating it, we obtain
\[w_{\epsilon2,0}(s)=w_{\epsilon2,0}(0)e^{2\lambda_\epsilon s}-\frac{1}{\mu_\epsilon+\omega_\epsilon
i}g_{\epsilon2,0}\left(e^{\lambda_\epsilon s}-e^{2\lambda_\epsilon
s}\right)-\]
\[-\frac{1}{\mu_\epsilon+3\omega_\epsilon
i}\overline{g}_{\epsilon0,2}\left(e^{\overline{\lambda}_\epsilon
s}-e^{2\lambda_\epsilon s}\right).
\]
By taking $s=-r$ we get
\[-e^{-2\lambda_\epsilon r}w_{\epsilon2,0}(0)+w_{\epsilon2,0}(-r)=-\frac{1}{\mu_\epsilon+\omega_\epsilon
i}g_{\epsilon2,0}\left(e^{-\lambda_\epsilon r}-e^{-2\lambda_\epsilon
r}\right)-\]
\[-\frac{1}{\mu_\epsilon+3\omega_\epsilon
i}\overline{g}_{\epsilon0,2}\left(e^{-\overline{\lambda}_\epsilon
r}-e^{-2\lambda_\epsilon r}\right),
\]
while the supplementary  condition is
\[
(A_{\epsilon}-2\lambda_\epsilon)
w_{\epsilon2,0}(0)+B_{\epsilon}w_{\epsilon2,0}(-r)=
g_{\epsilon2,0}+\overline{g}_{\epsilon0,2}-f_{\epsilon2,0}.
\]

\smallskip
The differential equation for $w_{\epsilon1,1}$ is
\[\frac{d w_{\epsilon1,1}(s)}{ds}=2\mu_\epsilon w_{\epsilon1,1}(s)+g_{\epsilon1,1}e^{\lambda_\epsilon
s}+ \overline{g}_{\epsilon1,1}e^{\overline{\lambda}_\epsilon s},
\]
from where,
\[w_{\epsilon1,1}(s)=w_{\epsilon1,1}(0)e^{2\mu_\epsilon s}-\frac{1}{\overline{\lambda}_\epsilon}g_{\epsilon1,1}\left(e^{\lambda_\epsilon
s}-e^{2\mu_\epsilon
s}\right)-\frac{1}{\lambda_\epsilon}\overline{g}_{\epsilon1,1}\left(e^{\overline{\lambda}_\epsilon
s}-e^{2\mu_\epsilon s}\right).
\]
We set again $s=-r,$
\[-e^{-2\mu_\epsilon r}w_{\epsilon1,1}(0)+w_{\epsilon1,1}(-r)=-\frac{1}{\overline{\lambda}_\epsilon}g_{\epsilon1,1}\left(e^{-\lambda_\epsilon
r}-e^{-2\mu_\epsilon
r}\right)-\frac{1}{\lambda_\epsilon}\overline{g}_{\epsilon1,1}\left(e^{-\overline{\lambda}_\epsilon
r}-e^{-2\mu_\epsilon r}\right)
\]
and by using also the condition for $w_{\epsilon1,1}:$
\[(A_{\epsilon}-2\mu_\epsilon)w_{\epsilon1,1}(0)+B_{\epsilon}w_{\epsilon1,1}(-r)
=g_{\epsilon1,1}+\overline{g}_{\epsilon1,1}-f_{\epsilon1,1},\] we
have the two equations for
$w_{\epsilon1,1}(-r),\,w_{\epsilon1,1}(0).$

For $w_{\epsilon0,2},$ it is enough to remark that
$w_{\epsilon0,2}=\overline{w_{\epsilon2,0}}.$

\end{document}